\documentclass[11pt]{article}

\usepackage[T1]{fontenc}
\usepackage{lmodern}
\usepackage{microtype}
\usepackage[margin = 1.28in]{geometry}
\usepackage{amsmath,amssymb,amsthm,mathtools}
\usepackage{enumitem}
\usepackage[colorlinks=true,linkcolor=blue,citecolor=blue,urlcolor=blue]{hyperref}

\newtheorem{theorem}{Theorem}[section]
\newtheorem{conjecture}[theorem]{Conjecture}
\newtheorem{lemma}[theorem]{Lemma}

\theoremstyle{remark}

\newcommand{\RR}{\mathbb{R}}
\newcommand{\aff}{\operatorname{aff}}

\newcommand{\congruent}{\mathrel{\cong}}

\newcommand{\base}{B}
\newcommand{\basept}{b}
\newcommand{\basedim}{{d_1}}
\newcommand{\target}{X}
\newcommand{\rcatcher}{C}
\newcommand{\catcherpt}{c}
\newcommand{\catcherdim}{{d_2}}
\newcommand{\catchersize}{{d_3}}
\newcommand{\catchersimplex}{S}
\newcommand{\simplexpt}{s}
\newcommand{\baseprod}{P}
\newcommand{\apexes}{A}
\newcommand{\apexpt}{a}

\title{A pyramid with a Ramsey base is Ramsey}

\author{Kenneth Moore\thanks{R\'enyi Institute, 1053 Budapest, Re\'altanoda u. 13-15, Hungary. Supported by ERC Advanced Grants 882971 ``GeoScape'' and ``ERCiD.'' Email:
{\tt moore.kenneth@renyi.hu}.}
}

\date{August 2026}

\begin{document}

\maketitle

\begin{abstract}
A finite subset $X$ of Euclidean space is called Ramsey if for any number of colours $k$ there exists a dimension $n$ such that whenever $\RR^n$ is $k$-coloured there exists a monochromatic congruent copy of $X$. The classification of Ramsey sets is one of the major unsolved problems in the field of Euclidean Ramsey theory. Towards this, Ivan, Leader and Walters recently asked whether adding a point to a Ramsey set outside of its affine hull necessarily produces another Ramsey set. In this note, we answer their question in the affirmative.
\end{abstract}

\section{Introduction}

A finite set $T$ in a Euclidean space is called \emph{Ramsey} if, for every integer $r\geq 1$, there is an integer $n$ such that every $r$-colouring of $\RR^n$ contains a monochromatic congruent copy of $T$. This notion was introduced by Erd\H{o}s, Graham, Montgomery, Rothschild, Spencer and Straus in three foundational papers \cite{EGMRSS1973,EGMRSS1975II,EGMRSS1975III}, where they established the field of \emph{Euclidean Ramsey theory}. 

One of the most important open problems in Euclidean Ramsey theory is to classify the Ramsey sets. Interestingly, it is known that all Ramsey sets must be spherical, meaning that they are contained in the surface of a sphere in some dimension. However, all of the known constructions of Ramsey sets have an even stronger property: subtransitivity. A \emph{transitive set} is a set whose group of symmetries acts transitively, and a subtransitive set is a subset of a transitive set. On this problem, we have

\begin{conjecture} 
\label{conj:main}
We give two competing statements.
\begin{itemize}
\item[\textbf{A.}] The Ramsey sets are the spherical sets.
\item[\textbf{B.}] The Ramsey sets are the subtransitive sets.
\end{itemize}
\end{conjecture}

Conjecture~\ref{conj:main}(A) was one of the main problems in the field for many years, originating from one of the three aforementioned papers of Erd\H{o}s et al. But a more recent work of Leader, Russell, and Walters \cite{LeaderRussellWalters2012} explains why Conjecture~\ref{conj:main}(B) may be much more natural. See also \cite{LeaderRussellWalters2011,Eberhard2013} for more on these conjectures.

In their recent paper, Ivan, Leader, and Walters \cite{IvanLeaderWalters2026} prove several results on prism and pyramid constructions that have symmetry properties. In particular, their Corollary~4 shows that a one-point extension of a transitive configuration produces a subtransitive configuration; under the corresponding solubility hypothesis, the resulting configuration is Ramsey. This gave a new method for constructing Ramsey configurations from subsoluble ones.

They then ask if the assumptions on the base can be replaced by the bare assumption that the base is Ramsey. In this note, we answer their question affirmatively with

\begin{theorem}[Conjecture 8 in \cite{IvanLeaderWalters2026}]
\label{thm:main}
Let $\base$ be a finite Ramsey set in a Euclidean space and let $z$ be a vector outside of the affine hull of $\base$. Then $\base\cup\{z\}$ is Ramsey.
\end{theorem}

This theorem would also follow from either of the competing conjectures in \ref{conj:main}, hence it does not support either conjecture. In \cite[Section 4]{IvanLeaderWalters2026} they note this fact, and also pose questions about other simple configurations that could produce evidence for one conjecture or the other.

The proof uses three classical results: the product theorem of Erd\H{o}s et al., according to which a Cartesian product of Ramsey configurations is Ramsey \cite{EGMRSS1973}, the theorem of Frankl and R\"odl that every nondegenerate simplex is Ramsey \cite{FranklRodl1990}, and a compactness consequence of the de Bruijn--Erd\H{o}s theorem. Apart from these inputs, the argument is elementary.

\section{Preliminaries}

For a set $S$ in a Euclidean space, we use $\operatorname{aff}(S)$ to denote the \emph{affine hull} of $S$, which is the set of all affine combinations of vectors in $S$. A \emph{congruent copy} of $S$ is another set $S'$ in a Euclidean space such that there exists a \emph{congruence}, which is a bijective map $\varphi:S \longrightarrow S'$ that preserves all pairwise distances. Note that if $S,S'\subseteq \RR^N$, $\varphi$ may be extended to an isometry of $\RR^N$. Given a second Euclidean set $T$, we employ the standard notation
$$
    S\xrightarrow{r} T
$$
if every $r$-colouring of $S$ contains a monochromatic congruent copy of $T$. Typically, this notation is used only when $T$ is a finite set. We will need two classical Euclidean Ramsey theorems, starting with

\begin{theorem}[Product theorem \cite{EGMRSS1973}]
\label{thm:product}
If $A$ and $B$ are Ramsey configurations, then their Cartesian product
$$
    A\times B=\{(a,b):a\in A,\ b\in B\}
$$
is Ramsey.
\end{theorem}

\begin{theorem}[Simplex theorem \cite{FranklRodl1990}]
\label{thm:simplex}
Every finite affinely independent Euclidean configuration is Ramsey.
\end{theorem}

We also use an implication of the standard hypergraph form of the de Bruijn--Erd\H{o}s theorem \cite{deBruijnErdos1951}. This is essentially Proposition~4 in \cite{EGMRSS1973}, which is also described at the top of Section~2 in \cite{LeaderRussellWalters2012}. We state it as

\begin{lemma}\label{lem:witness} For any finite configuration $T$, if $\RR^n\xrightarrow{r} T,$ then there exists a finite set $A\subseteq\RR^n$ such that $A\xrightarrow{r} T.$
\end{lemma}

\section{Proof of Theorem~\ref{thm:main}}

\label{section:proof}

Denote the target set by $\target=\base \cup\{z\}.$ We prove by induction on $r\geq 1$ that $\target$ is Ramsey for $r$ colours. The case when $r=1$ is not difficult.

\vspace{3mm}

We first give a brief intuitive explanation of the argument. We are going to find a finite configuration $\rcatcher$ such that every $(r-1)$-colouring of $\rcatcher$ contains a monochromatic copy of $\target$. We next create a set $\baseprod$ as a product of the base $\base$ with a simplex, which is a Ramsey set. The simplex is chosen carefully, to make every copy of $\base$ `aim' at a different point of a particular copy of $\rcatcher$, such that adding that point to $\base$ would create a copy of $\target$. 

We choose the ambient dimension sufficiently large so that we can find a monochromatic (say, red) copy of $\baseprod$. Now, if any point $\catcherpt$ of that nearby copy of $\rcatcher$ is red, there is a subset of $\baseprod$ congruent to $\base$ which combines with $\catcherpt$ to create a red $\target$, and we are done. Otherwise, no point of the copy of $\rcatcher$ is red, and in particular, only $r-1$ colours are used. Thus we can find a monochromatic copy of $\target$ by induction. 

\begin{figure}[ht]
\centering
\includegraphics[scale=0.55]{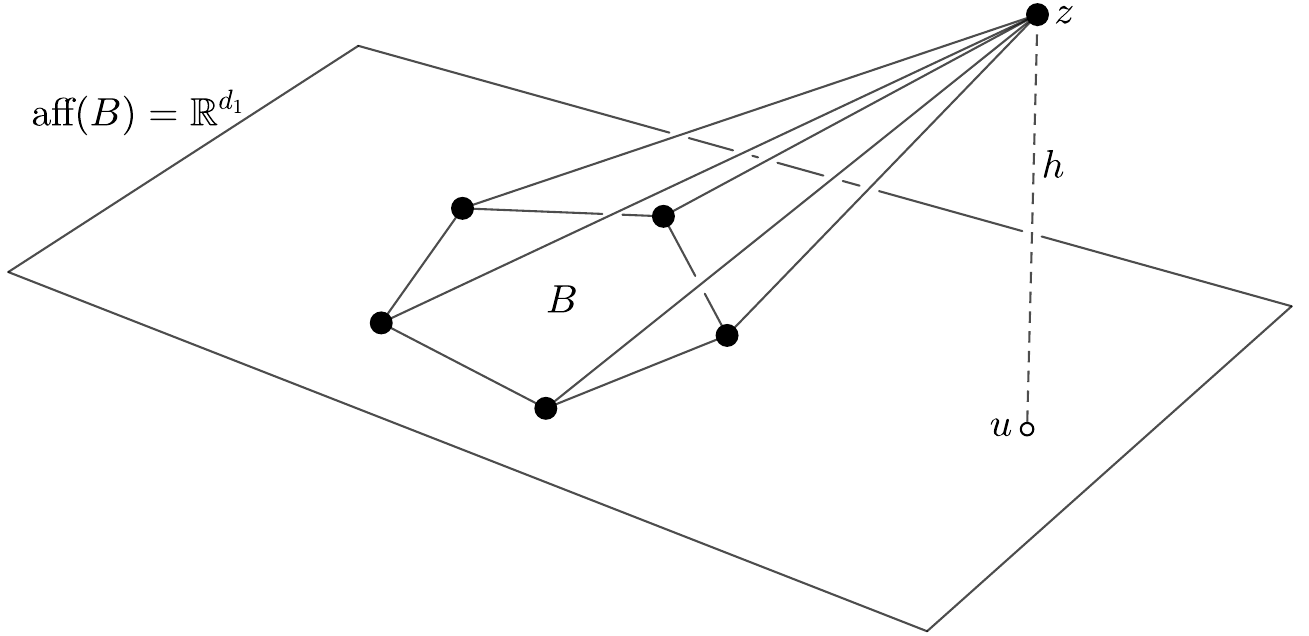}
\caption{The pyramid set up $Y=\base \cup\{z\}$}
\label{fig:pyramid}
\end{figure}

We now proceed with the proof proper. Fix $r\geq 2$ and assume that $\target$ is known to be Ramsey for $r-1$ colours; so there is a dimension $\catcherdim$ such that
$$\RR^\catcherdim \xrightarrow{r-1} \target.$$
By Lemma~\ref{lem:witness}, there is a finite configuration
$$\rcatcher=\{\catcherpt_1,\dots,\catcherpt_\catchersize\}\subseteq \RR^\catcherdim$$
for some integer $\catchersize$ satisfying
\begin{equation}\label{eq:A-witness}
    \rcatcher\xrightarrow{r-1} \target.
\end{equation}

After applying an isometry if necessary, we may assume $\aff(\base)$ is the Euclidean space $\RR^\basedim$, and that we can decompose $z$ into its orthogonal projection $u\in \RR^\basedim$ and its perpendicular component $h\in\RR$, with $h>0$. So $z=(u,h)\in \RR^{\basedim+1},$ and $h>0$ (as in Figure~\ref{fig:pyramid}). Consequently, for all $\basept\in \base$,
\begin{equation}\label{eq:Y-distances}
    \|(\basept,0)-z\|^2=\|\basept-u\|^2+h^2.
\end{equation}

Let $e_1,\dots,e_\catchersize$ be an orthonormal basis of $\RR^\catchersize$. Define $$\simplexpt_i=(\catcherpt_i,h e_i)\in \RR^{\catcherdim+\catchersize}\qquad\text{for}\qquad 1\leq i\leq \catchersize,$$ and let
\begin{equation*}
    \catchersimplex =\{\simplexpt_1,\dots,\simplexpt_\catchersize\}.
\end{equation*}
$\catchersimplex$ is clearly affinely independent, since the $he_i$-coordinates are formed by linearly independent vectors. Therefore, Theorem~\ref{thm:simplex} implies that $\catchersimplex$ is Ramsey.

Theorem~\ref{thm:product} now shows that the product set
\begin{equation*}
    \baseprod:=\base\times \catchersimplex
      = \{(\basept,\catcherpt_i,h e_i):\basept\in \base,\ 1\leq i\leq \catchersize\}
      \subseteq \RR^{\basedim + \catcherdim + \catchersize}
\end{equation*}
is Ramsey. For each $i$, introduce the corresponding apex point
\begin{equation*}
    \apexpt_i=(u,\catcherpt_i,0)\in \RR^{\basedim + \catcherdim + \catchersize},
\end{equation*}
and let
$$ \apexes=\{\apexpt_1,\dots,\apexpt_\catchersize\}.
$$

Since $u$ is constant, the configuration $\apexes$ is congruent to $\rcatcher$. For each fixed $i$, we define the fibre
$$\baseprod_i := \base\times\{\simplexpt_i\} =\{(\basept,\catcherpt_i,h e_i):\basept\in \base\}
$$
which is congruent to $\base$. Further, for every $\basept\in \base$,
\begin{equation}
\begin{split} 
    \|(\basept,\catcherpt_i,h e_i)-\apexpt_i\|^2
      &=\|\basept-u\|^2+\|h e_i\|^2 \\
      &=\|\basept-u\|^2+h^2.
\end{split} 
\label{eq:fibre-apex-distance}
\end{equation}
Comparing \eqref{eq:fibre-apex-distance} with \eqref{eq:Y-distances}, we obtain
\begin{equation}\label{eq:fibre-Y}
    \baseprod_i\cup\{\apexpt_i\}\congruent \target
    \qquad (1\leq i\leq \catchersize ).
\end{equation}

Because $\baseprod$ is Ramsey, there is a dimension $n_0$ such that $ \RR^{n_0}\xrightarrow{r}\baseprod.$ Let $$N=\max \{n_0,\basedim + \catcherdim + \catchersize\}$$ and consider $P$ and $C$ to be subsets of $\RR^N$. For an arbitrary $r$-colouring of $\RR^N$, there is a monochromatic copy $P'\subset\RR^N$ of $P$. Call its colour red, and let
$$
    \varphi:\baseprod \longrightarrow \baseprod'
$$
be the corresponding congruence. We extend $\varphi$ to an isometry $\widetilde\varphi$ on $\RR^N$, and write
$$
   \apexpt_i'=\widetilde\varphi(\apexpt_i).
$$
Then every point of every fibre $\widetilde\varphi(\baseprod_i)$ is red.

\newpage 

There are two cases. If $\apexpt_i'$ is red for some $i$, then
$$
    \widetilde\varphi(\baseprod_i\cup\{\apexpt_i\})
$$
is a red copy of $\target$ by \eqref{eq:fibre-Y}. Otherwise none of
$\apexpt_1',\dots,\apexpt_\catchersize'$ is red. These points therefore use at most the remaining
$r-1$ colours. Since $\widetilde\varphi$ is an isometry, they form a congruent copy of $\rcatcher$; by \eqref{eq:A-witness}, they contain a monochromatic copy of $\target$. This completes the proof of Theorem~\ref{thm:main}.

\section*{Acknowledgements}

ChatGPT 5.6 (OpenAI) was used during the brainstorming stages of this article.

\bibliographystyle{abbrv}
\bibliography{bibliography}

\end{document}